\documentclass[letterpaper,10pt]{amsart}
\usepackage{amsmath,amsthm,amstext}
\usepackage{amsfonts,amssymb}
\usepackage[mathscr]{eucal}
\usepackage[]{hyperref}
\usepackage{setspace}

\newtheorem{theorem}{Theorem}[section]

\theoremstyle{definition}
\newtheorem{remark}[theorem]{Remark}
\newtheorem{definition}[theorem]{Definition}

\numberwithin{equation}{section}

\renewcommand{\l}{\lambda}

\newcommand{\RR}{\ensuremath{\mathbb{R}}}
\newcommand{\CC}{\ensuremath{\mathbb{C}}}
\newcommand{\sph}{\ensuremath{\mathbb{S}}}
\newcommand{\prtl}{\ensuremath{\partial}}
\DeclareMathOperator{\supp}{supp}

\newcommand{\g}{\ensuremath{\mathrm{g}}}
\newcommand{\h}{\ensuremath{\mathrm{h}}}
\newcommand{\To}{\longrightarrow}
\newcommand{\upd}{\ensuremath{\mathrm{d}}}
\newcommand{\defeq}{\ensuremath{\stackrel{\mathrm{def}}{=}}}

\DeclareMathOperator{\dist}{dist}

\title[Strichartz estimates on polygonal domains]{Strichartz estimates
  for the Schr\"odinger equation on polygonal domains}

\thanks{MDB was supported by the NSF grant DMS-0801211.  GAF was partially supported by NSF grant DMS-0636646.  JLM was supported in part by an NSF Postdoctoral Fellowship and a Hausdorff Center Postdoc at the University of Bonn.}

\author[M. D. Blair]{Matthew D. Blair}
\address{Department of Mathematics and Statistics, University of New Mexico, Albuquerque, NM 87131, USA}
\email{blair@math.unm.edu}

\author[G. A. Ford]{G. Austin Ford}
\address{Department of Mathematics, Northwestern University, Evanston, IL 60208, USA}
\email{aford@math.northwestern.edu}

\author[S. Herr]{Sebastian Herr}
\address{Mathematisches Institut, Universit\"at Bonn, 53115 Bonn, Germany}
\email{herr@math.uni-bonn.de}

\author[J. L. Marzuola]{Jeremy L. Marzuola}
\address{Applied Mathematics Department, Columbia University, New York, NY 10027, USA}
\email{jm3058@columbia.edu}

\begin{document}

\maketitle

\begin{abstract}
We prove Strichartz estimates with a loss of derivatives for the Schr\"odinger equation on polygonal domains with either Dirichlet or Neumann homogeneous boundary conditions.  Using a standard doubling procedure, estimates the on polygon follow from those on Euclidean surfaces with conical singularities.  We develop a Littlewood-Paley squarefunction estimate with respect to the spectrum of the Laplacian on these spaces.  This allows us to reduce matters to proving estimates at each frequency scale.  The problem can be localized in space provided the time intervals are sufficiently small.  Strichartz estimates then follow from a result of the second author regarding the Schr\"odinger equation on the Euclidean cone.
\end{abstract}

\section{Introduction}\label{sec:Intro}

Let $\Omega \subset \RR^2$ be a compact polygonal domain in the
plane, that is, a compact, connected region in $\RR^2$ whose
boundary is piecewise linear. Suppose $u(t,x): [-T,T] \times \Omega \To
\CC$ is a solution to the initial value problem for the Schr\"{o}dinger equation on $\Omega$,

\begin{equation}\label{eq:Schro IVP}
\left\{
\begin{split}
\left(D_t + \Delta\right) u(t,x) &= 0 \\
u(0,x) &= f(x)  ,
\end{split}
\right.
\end{equation}
satisfying either Dirichlet or Neumann homogeneous boundary conditions,
\begin{equation}\label{eq:bdry conditions}
u\big|_{[-T,T] \times \prtl \Omega} = 0 \qquad \text{ or } \qquad \prtl_n u\big|_{[-T,T] \times \prtl \Omega} = 0  .
\end{equation}
Here, $\prtl_n$ denotes the normal derivative along the boundary,
$D_t = \frac{1}{i} \, \prtl_t$, and $\Delta = - \prtl_{x_1}^2 - \prtl_{x_2}^2$ is the nonnegative Laplace operator.

In this note, we are interested in Strichartz estimates for solutions to the aforementioned Schr\"odinger IBVP \eqref{eq:Schro IVP}-\eqref{eq:bdry conditions}; these are a family of
space-time integrability bounds of the form
\begin{equation}\label{eq:Strichartz}
\|u\|_{L^p([-T,T]; L^q(\Omega))} \leq C_T \left\|f\right\|_{H^s(\Omega)}
\end{equation}
with $p > 2$ and $\frac{2}{p} + \frac{2}{q} = 1$. In this estimate, we take the space $H^s(\Omega)$ to be the $L^2$-based Sobolev space of order $s$ defined with respect to the spectral resolution of either the Dirichlet or Neumann Laplacian.  More precisely, this self-adjoint operator possesses a sequence of eigenfunctions forming a basis for $L^2(\Omega)$.  We write the eigenfunction and eigenvalue pairs as $\Delta \varphi_j = \l^2_j \varphi_j$, where $\l_j$ denotes the frequency of vibration.  The Sobolev space of order $s$ can then be defined as the image of $L^2(\Omega)$ under $(1+\Delta)^{-s}$ with norm
\begin{equation}\label{eq:Sobolev defn}
\left\|f\right\|_{H^s(\Omega)}^2 = \sum_{j=1}^\infty \left( 1 + \l_j^2 \right)^s \left| \left\langle f, \varphi_j \right\rangle \right|^2  .
\end{equation}
Here, $\langle \cdot, \cdot \rangle$ denotes the $L^2$ inner product.

Strichartz estimates are well-established when the domain $\Omega$ is replaced by Euclidean space.  In that case, one can take $s=0$ in \eqref{eq:Strichartz}, and by scaling considerations, this is the optimal order for the Sobolev space; see for example Strichartz~\cite{strich77}, Ginibre and Velo~\cite{ginvelo85}, Keel and Tao~\cite{keeltao98}, and references therein.  When $\Omega$ is a compact domain or manifold, much less is known about the validity and optimality of these estimates.  The finite volume of the manifold and the presence of trapped geodesics appear to limit the extent to which dispersion can occur. In addition, the imposition of
boundary conditions complicate many of the known techniques for proving Strichartz estimates.  Nonetheless, estimates on general compact domains with smooth boundary have been shown by Anton \cite{antonint} and Blair-Smith-Sogge \cite{BSSpams}.  Both of
these works build on the approach for compact manifolds of Burq-G\'erard-Tzvetkov~\cite{burq1}.

In the present work, we prove the following
\begin{theorem}\label{thm:Strichartz on polygonal domains}
Let $\Omega$ be a compact polygonal domain in $\mathbb{R}^2$, and let $\Delta$ denote either the Dirichlet or Neumann Laplacian on $\Omega$.  Then for any solution $u = \exp(-it\Delta) \, f$ to the Schr\"odinger IBVP \eqref{eq:Schro IVP}-\eqref{eq:bdry conditions} with $f$ in $H^{\frac{1}{p}}(\Omega)$, the Strichartz estimates
\begin{equation}\label{eq:Strichartz on polygonal domains}
\|u\|_{L^p([-T,T];L^q(\Omega))} \leq C_T \left\| f \right\|_{H^{\frac{1}{p}}(\Omega)}
\end{equation}
hold provided $p > 2$, $q \geq 2$, and $\frac{2}{p} + \frac{2}{q} = 1$.
\end{theorem}

\begin{remark}
In this work, the Neumann Laplacian is taken to be the Friedrichs extension of the Laplace operator acting on smooth functions which vanish in a neighborhood of the vertices and whose normal derivative is zero on the rest of the boundary.  In this sense, our Neumann Laplacian imposes Dirichlet conditions at the vertices and Neumann conditions elsewhere.
The Dirichlet Laplacian is taken to be the typical Friedrichs extension of the Laplace operator acting on smooth functions which are compactly supported in the interior of $\Omega$.
\end{remark}

\begin{remark}
We note that our estimates have a loss of $s=\frac{1}{p}$ derivatives as in \cite{burq1}, which we believe is an artifact of our methods.  Given specific geometries, there are results showing that such a loss is not sharp.  For instance, when $\Omega$ is replaced by a flat rational torus, the estimate \eqref{eq:Strichartz} with $p=q=4$ holds for any $s>0$, as was shown by Bourgain~\cite{Bourgain93}; see also \cite{Bourgain07} for results in the case of irrational tori.  However, we also point out that in certain geometries a loss of derivatives is expected due to the existence of gliding rays, as shown by Ivanovici \cite{iv1}.
\end{remark}

\begin{remark}
Using a now standard application of the Christ-Kiselev lemma \cite{ChKi}, we can conclude that for a solution $u$ to the inhomogeneous Schr\"odinger IBVP
\begin{equation}\label{eq:inhom Schro IVP}
\left\{
\begin{split}
\left(D_t + \Delta\right) u(t,x) &= F(t,x) \\
u(0,x) &= f(x)
\end{split}
\right.
\end{equation}
satisfying either Dirichlet or Neumann homogeneous boundary
conditions, the estimate
\begin{equation}\label{eq:inhom Strichartz on polygonal domains}
\|u\|_{L^{p_1}([-T,T]; L^{q_1}(\Omega))} \leq C_T \left( \left\|f\right\|_{H^{\frac{1}{p_1}}(\Omega)} + \left\| F \right\|_{L^{p_2'}([-T,T];W^{\frac{1}{p_1}+\frac{1}{p_2},q_2'}(\Omega))} \right)  ,
\end{equation}
holds for $\frac{2}{p_j} + \frac{2}{q_j} = 1$ for $j = 1,2$.  Here, $(\cdot)'$ denotes the dual exponent, e.g.\ $\frac{1}{p_1} + \frac{1}{p_1'} = 1$.
\end{remark}

We prove Theorem \ref{thm:Strichartz on polygonal domains} by utilizing a doubling procedure to reduce estimates on the polygonal domain $\Omega$ to estimates on a Euclidean surface with conical singularities.  A Euclidean surface with conical singularities (ESCS) is, loosely speaking, a Riemannian surface $(X,\g)$ locally modeled on either Euclidean space or the flat Euclidean cone; for a precise definition, see Section \ref{sec:ESCSs}.  As will be outlined below, any compact planar polygonal domain $\Omega$ can be doubled across its boundary to produce a compact ESCS.  In this procedure, a vertex of $\Omega$ of angle $\alpha$ gives rise to a conic point of $X$ with cone angle $2\alpha$.  Taking the Laplacian on $X$ to be the Friedrichs extension of the Laplacian on $\mathcal{C}^\infty_c(X_0)$, where $X_0$ is $X$ less the singular points, and the Sobolev spaces $H^s(X)$ as in \eqref{eq:Sobolev defn}, Theorem \ref{thm:Strichartz on polygonal domains} will follow from the following

\begin{theorem}\label{thm:Strichartz on ESCSs}
Let $X$ be a compact ESCS, and let $\Delta_\g$ be the Friedrichs extension of $\Delta_\g \Big\vert_{\mathcal{C}^\infty_c(X_0)}$.  Then for any solution $u = \exp\!\left(-it\Delta_\g\right) f$ to the Schr\"odinger IVP on $X$ with initial data $f$ in $H^{\frac{1}{p}}(X)$, the Strichartz estimates
\begin{equation}\label{eq:Strichartz on ESCSs}
\|u\|_{L^p([-T,T];L^q(X))} \leq C_T \left\| f \right\|_{H^{\frac{1}{p}}(X)}
\end{equation}
hold provided $p > 2$, $q \geq 2$, and $\frac{2}{p} + \frac{2}{q} = 1$.
\end{theorem}

The method here is to develop a local parametrix for the operator at frequency-dependent scales using a Littlewood-Paley decomposition.  Since an ESCS locally looks like either the plane or the Euclidean cone, estimates will follow from a result of the second author~\cite{ford}, which develops Strichartz estimates on the latter.  However, since propagation speed is proportional to frequency, the error in the parametrix is only bounded over time intervals of size inversely proportional to the frequency scale, cp.\ \cite{burq1}.  The loss of $\frac{1}{p}$ derivatives relative to estimates on the plane thus results from decomposing the time interval $[-T,T]$ into smaller frequency-dependent time intervals over which the error is bounded.

\subsection*{Acknowledgement.} The authors are grateful to Luc Hillairet for valuable discussions regarding Euclidean surfaces with conical singularities.

\section{Euclidean surfaces with conical singularities}
\label{sec:ESCSs}

In this section, we review the definition and properties of Euclidean surfaces $X$ with conical singularities.  For more details, we refer the reader to \cite{C79}, where we believe manifolds of this type to have been first introduced, and \cite{Hil} and \cite{HHM}, where properties and applications of such surfaces are explored.

We begin by establishing the notation $C(\mathbb{S}^1_\rho) \defeq \mathbb{R}_+ \times \left(\mathbb{R} \big/ 2\pi\rho \mathbb{Z}\right)$ for the flat Euclidean cone of radius $\rho > 0$ equipped with the metric $\h(r,\theta) = \upd r^2 + r^2 \, \upd\theta^2$.  With this in mind, we have the following

\begin{definition}
A \emph{Euclidean surface with conical singularities} (ESCS) is a topological space $X$ possessing a decomposition $X = X_0 \sqcup P$ for a finite set of singular points $P \subsetneq X$ such that
\begin{enumerate}
	\item  $X_0$ is an open, smooth two-dimensional Riemannian manifold with a locally Euclidean metric $\g$, and
	\item  each singular point $p_j$ of $P$ has a neighborhood $U_j$ such that $U_j \setminus \left\{p_j\right\}$ is isometric to a neighborhood of the tip of a flat Euclidean cone $C(\mathbb{S}^1_{\rho_j})$ with $p_j$ mapped to the cone tip.
\end{enumerate}
\end{definition}

We stress that the analysis required to prove Theorem~\ref{thm:Strichartz on ESCSs} all occurs on the Riemannian manifold $X_0$.  As remarked previously, we take the Laplacian $\Delta_\g$ on $X$ to be the Friedrichs extension of $\Delta_\g \Big\vert_{\mathcal{C}^\infty_c(X_0)}$.  This is a nonnegative, self-adjoint operator on $L^2(X)$ with discrete spectrum tending to infinity, as can be seen from the Rellich-type theorem of \cite[Theorem 3.4]{CT82}.  We can thus take the Sobolev spaces $H^s(X)$ to be the images of $L^2(X)$ under $\left(1 + \Delta_\g\right)^{-s}$ with norm defined similarly to that in \eqref{eq:Sobolev defn}.

We now discuss how any compact polygonal domain $\Omega$ in $\mathbb{R}^2$, possibly with polygonal holes, gives rise to an ESCS $X$ equipped with a flat metric $\g$.  Begin with two copies $\Omega$ and $\sigma \Omega$ of the polygonal domain, where $\sigma$ is a reflection of the plane. The double $X$ is obtained by taking the formal union $\Omega \cup \sigma \Omega$, where two corresponding sides are identified pointwise.  Taking polar coordinates near each vertex of the polygon, it can be seen that the flat metric $\g$ extends smoothly across the sides.  In particular, a vertex in $\Omega$  of angle $\alpha$ gives rise to a conic point of $X$ locally isometric to $C(\mathbb{S}^1_\rho)$ with $\rho = \frac{\alpha}{\pi}$.

The reflection $\sigma$ of $\Omega$ gives rise to an involution of $X$ commuting with the Laplace operator.  This operator $\Delta_\g$ thus decomposes into two operators acting on functions which are odd or even with respect to $\sigma$, and these operators are then equivalent to the Laplace operator on $\Omega$ with Dirichlet or Neumann boundary conditions respectively.  In particular, for any eigenfunction $\varphi_j$ of the Dirichlet, resp.\ Neumann, Laplace operator on $\Omega$, we can construct an eigenfunction of the Laplace operator on $X$ by taking $\varphi_j$ in $\Omega$ and $-\varphi_j \circ \sigma$, resp. $\varphi_j \circ \sigma$, in $\sigma \Omega$.  As a consequence, the Schr\"odinger flow over $X$ can be seen to extend that for $\Omega$, and hence the Strichartz estimates in Theorem~\ref{thm:Strichartz on polygonal domains} follow from those in Theorem~\ref{thm:Strichartz on ESCSs} as claimed.

\section{Strichartz estimates}

In this section we prove Theorem~\ref{thm:Strichartz on ESCSs}.  We start with a Littlewood-Paley decomposition of our solution $u$ in the spatial frequency domain.  Namely, choose a nonnegative bump function $\beta$ in $\mathcal{C}^\infty_c(\mathbb{R})$ supported in $\left(\frac{1}{4},4\right)$ and satisfying $\sum_{k \geq 1} \beta\!\left(2^{-k} \, \zeta\right) = 1$ for $\zeta \geq 1$.  Taking $\beta_k(\zeta) \defeq \beta\!\left(2^{-k} \, \zeta\right)$ for $k \geq 1$ and $\beta_0(\zeta) \defeq 1 - \sum_{k \geq 1} \beta_k(\zeta)$, we define the frequency localization $u_k$ of $u$ in the spatial variable by
\begin{equation}\label{eq:L-P localization}
u_k \defeq \beta_k\!\left(\sqrt{\Delta_\g}\right) u  ,
\end{equation}
where the operator $\beta_k\!\left(\sqrt{\Delta_g}\right)$ is defined using the functional calculus with respect to $\Delta_g$.  Hence, $u = \sum_{k \geq 0} u_k$, and in particular, $u_0$ is localized to frequencies smaller than $1$.

With this decomposition, we have the following squarefunction estimate for elements $a$ of $L^q(X)$,
\begin{equation}\label{eq:L-P estimate}
\left\| \left(\sum_{k \geq 0} \left|\beta_k\!\left(\sqrt{\Delta_\g} \right) a \right|^2\right)^{\frac{1}{2}} \right\|_{L^q(X)} \approx \left\|a\right\|_{L^q(X)}  ,
\end{equation}
with implicit constants depending only on $q$.  Delaying the proof of \eqref{eq:L-P estimate} to Section \ref{sec:L-P squarefn estimate}, we have by Minkowski's inequality that
\begin{equation}\label{eq:Minkowski}
\|u\|_{L^p([-T,T]; L^q(X))} \lesssim \left( \sum_{k \geq 0} \left\|u_k\right\|_{ L^p([-T,T]; L^q(X))}^2 \right)^{\frac{1}{2}}
\end{equation}
since we are under the assumption that $p,q \geq 2$.  We now claim that for each $k \geq 0$,
\begin{equation}\label{eq:dyadic Strichartz}
\left\|u_k\right\|_{L^p([-T,T]; L^q(X))} \lesssim 2^{\frac{k}{p}} \left\|u_k(0,\cdot)\right\|_{L^2(X)} .
\end{equation}

Assuming this for the moment, we have by orthogonality and the localization of $\beta$ that
\begin{equation}
\begin{aligned}
2^{\frac{2k}{p}} \left\|u_k(0,\cdot)\right\|_{L^2(X)}^2 &= 2^{\frac{2k}{p}} \sum_{j=1}^\infty \beta_k(\l_j)^2 \left| \left\langle u(0,\cdot),\varphi_j \right\rangle \right|^2 \\
& \lesssim \sum_{j=1}^\infty \left(1 + \l_j^2\right)^{1/p} \beta_k(\l_j)^2 \left| \left\langle u(0,\cdot),\varphi_j \right\rangle \right|^2.
\end{aligned}
\end{equation}
We now sum this expression over $k$; after exchanging the order of summation in $k$ and $j$, we obtain
\begin{equation}
\sum_{k \geq 0} 2^{\frac{2k}{p}} \left\|u_k(0,\cdot)\right\|_{L^2(X)}^2 \lesssim \left\|u(0,\cdot)\right\|_{H^{\frac{1}{p}}(X)}^2.
\end{equation}
Combining this with~\eqref{eq:Minkowski}, we have reduced the proof of Theorem \ref{thm:Strichartz on ESCSs} to showing the claim~\eqref{eq:dyadic Strichartz}.

We now observe that~\eqref{eq:dyadic Strichartz} follows from
\begin{equation}\label{eq:semiclass}
\|u_k\|_{ L^p ([0,2^{-k}]; L^q(X))} \lesssim \|u_k(0,\cdot)\|_{L^2(X)}.
\end{equation}
Indeed, if this estimate holds, then time translation and mass conservation imply the same estimate holds with the time interval $[0, 2^{-k}]$ replaced by $[2^{-k}m, 2^{-k}(m+1)] $.  Taking a sum over all such dyadic intervals in $[-T,T]$ then yields~\eqref{eq:dyadic Strichartz}.

Next, we localize our solution in space using a finite partition of unity $\sum_\ell \psi_\ell \equiv 1$ on $X$ such that $\supp(\psi_\ell)$ is contained in a neighborhood $U_\ell$ isometric to either an open subset of the plane $\RR^2$ or a neighborhood of the tip of a Euclidean cone $C(\mathbb{S}^1_\rho)$.  It now suffices to see that if $\psi$ is an element of this partition and $U$ denotes the corresponding open set in $\RR^2$ or $C(\mathbb{S}^1_\rho)$, then
\begin{equation}
\left\|\psi \, u_k\right\|_{L^p([0,2^{-k}];L^q(U))} \lesssim \left\|u_k(0,\cdot)\right\|_{L^2(U)}.
\end{equation}
Here and in the remainder of the section, $L^q(U)$ is taken to mean the space of functions on $U$ which are $q$-integrable with respect to the Riemannian measure over $\RR^2$ or $C(\mathbb{S}^1_\rho)$, depending on where $U$ lies.

Observe that $\psi \, u_k$ solves the equation
\begin{equation}
\left(D_t + \Delta_\g\right) \left(\psi\, u_k \right) = \left[\Delta_\g, \psi\right] u_k
\end{equation}
over $\RR^2$ or $C(\mathbb{S}^1_\rho)$.  Letting $\mathbf{S}(t)$
denote the Schr\"odinger propagator either on Euclidean space or the
Euclidean cone, depending on which space $U$ lives in, we have for $t \geq 0$ that
\begin{equation}
\psi \, u_k(t,\cdot) = \mathbf{S}(t) \big(\psi\, u_k(0,\cdot)\big) + \int_{0}^{2^{-k}} \mathbf{1}_{ \{t>s\} }(s) \, \mathbf{S}(t-s)\!\left(\left[\Delta_\g,\psi\right] u_k(s,\cdot) \right) \, \upd s  .
\end{equation}
Here, $\mathbf{1}_{ \{t>s\} } (s)$ is the indicator of the set $t>s>0$. By Minkowski's inequality,
\begin{multline}
\left\|\psi \, u_k \right\|_{L^p([0,2^{-k}];L^q(U))} \lesssim \left\| \, \mathbf{S}(\cdot)  \big(\psi \, u_k(0,\cdot) \big) \right\|_{L^p([0,2^{-k}];L^q(U))} \\
\mbox{} + \int_{0}^{2^{-k}} \left\| \, \mathbf{S}(\cdot-s) \! \left(\left[\Delta_\g,\psi\right] u_k(s,\cdot) \right) \right\|_{L^p([0,2^{-k}];L^q(U))} \, \upd s  .
\end{multline}

We now apply known Strichartz estimates on $\mathbf{S}(t)$.  When $U$ is a subset of the plane, the estimates on the propagator are well-known and contained in the references listed in the introduction.  When $U$ is a subset of the flat Euclidean cone $C(\sph^1_\rho)$, the estimates are due to the following result of the second author.
\begin{theorem}[Theorem 5.1 of \cite{ford}]
Suppose $p>2$ and $q \geq 2$ satisfy $\frac{2}{p} + \frac{2}{q} = 1$.  Then the Schr\"odinger solution operator $\mathbf{S}(t) =  \exp(-it\Delta_\h)$ on $C(\sph^1_\rho)$ satisfies the Strichartz estimates
\begin{equation}
\|\mathbf{S}(t) f\|_{L^p(\RR;L^q(C(\sph^1_\rho)))} \lesssim \|f\|_{L^2(C(\sph^1_\rho))}  .
\end{equation}
\end{theorem}

\noindent We now conclude that
\begin{multline}\label{eq:local Strichartz}
\left\|\psi \, u_k\right\|_{L^p([0,2^{-k}];L^q(U))} \lesssim \left\|\psi \, u_k(0,\cdot)\right\|_{L^2(U)} \\
\mbox{} + \int_{0}^{2^{-k}} \left\| \left(\left[\Delta_\g,\psi \right] u_k(s,\cdot) \right) \right\|_{L^2(U)} \, \upd s  .
\end{multline}
The estimates~\eqref{eq:semiclass} will then follow provided
\begin{equation}\label{eq:commutator}
\begin{aligned}
2^{-k} \left\| \left[\Delta_\g, \psi \right] u_k \right\|_{L^\infty L^2(X)} &\lesssim 2^{-k} \left\|\nabla_\g u_k \right\|_{L^\infty L^2(X)} + 2^{-k} \left\|u_k \right\|_{L^\infty L^2(X)} \\
&\lesssim \left\|u_k(0,\cdot)\right\|_{L^2(X)}  ,
\end{aligned}
\end{equation}
with $\nabla_\g$ denoting the Riemannian gradient.  The first inequality here follows by a simple computation of the commutator. For the second inequality, first observe that the term $2^{-k}\|u_k\|_{L^\infty L^2}$ is easily controlled by mass conservation.  We then claim that the bound on the gradient term follows from
\begin{equation}\label{eq:Lapl bound time}
2^{-2k} \left\|\Delta_\g u_k \right\|_{L^\infty L^2(X)} \lesssim \left\|u_k(0,\cdot)\right\|_{L^2(X)}  .
\end{equation}
Indeed, if this holds we have that
\begin{equation}
\begin{aligned}
2^{-2k} \left\|\nabla_\g u_k\right\|_{L^\infty L^2(X)}^2 &\lesssim \sup_t 2^{-2k} \left\langle \Delta_\g u_k(t,\cdot), u_k(t,\cdot) \right\rangle \\
&\lesssim \sup_t 2^{-2k} \left\|\Delta_\g u_k(t,\cdot)\right\|_{L^2(X)} \left\|u_k(t,\cdot)\right\|_{L^2(X)} \\
&\lesssim \left\|u_k(0,\cdot)\right\|^2_{L^2(X)}.
\end{aligned}
\end{equation}
We next observe that since the Schr\"odinger propagator $\mathbf{S}(t) = \exp(-it\Delta_\g)$ commutes with $\Delta_\g$, mass conservation implies that~\eqref{eq:Lapl bound time} further reduces to showing the bound
\begin{equation}\label{eq:Lapl bound}
  2^{-2k} \left\|\Delta_\g \, \beta_k\!\left(\sqrt{\Delta_\g}\right) a \right\|_{L^2(X)} \lesssim \left\|a\right\|_{L^2(X)}
\end{equation}
for elements $a$ of $L^2(X)$.

Finally, to see that~\eqref{eq:Lapl bound} holds, define $\Psi(t)$ to be the
Schwartz class function satisfying
\begin{equation}
\zeta^2 \, \beta(\zeta) = \int_{-\infty}^{\infty} e^{it\zeta} \, \Psi(t) \, \upd t  .
\end{equation}
This implies
\begin{equation}
2^{-2k} \, \Delta_\g \, \beta_k\!\left(\sqrt{\Delta_\g}\right) = 2^k \int_{-\infty}^{\infty} e^{it\sqrt{\Delta_\g}} \, \Psi\!\left(2^{k}t\right) \, \upd t  .
\end{equation}
We now use that $\exp\!\left(it\sqrt{\Delta_\g}\right)$ is an isometry on
$L^2(X)$ to obtain
\begin{equation}
\begin{aligned}
\left\|2^{-2k} \, \Delta_\g \, \beta_k\!\left(\sqrt{\Delta_\g}\right) a \right\|_{L^2(X)} &\leq 2^k \int_{-\infty}^{\infty} \left\|e^{it\sqrt{\Delta_\g}}a\right\|_{L^2(X)} \Psi\!\left(2^{k}t\right) \upd t \\
&\lesssim \left\|a\right\|_{L^2(X)}  ,
\end{aligned}
\end{equation}
showing \eqref{eq:Lapl bound} and thus, moving backwards through the reductions, the claim \eqref{eq:dyadic Strichartz}.

\section{The Littlewood-Paley squarefunction estimate}
\label{sec:L-P squarefn estimate}

In this section, we prove the Littewood-Paley squarefunction estimate~\eqref{eq:L-P estimate} for Euclidean surfaces with conical singularities, which is the last remaining piece of the proof of Theorem \ref{thm:Strichartz on ESCSs}.  As we shall see, the estimate is actually valid for any exponent $1 < q <\infty$.  If $X_0$ were compact, the estimate in Seeger-Sogge~\cite[Lemma 2.3]{seegersogge} would suffice for our purpose.  Extra care must be taken in our case, however, as $X_0$ is an incomplete manifold.  Thus, we take advantage of a spectral multiplier theorem that allows us to employ a classical argument appearing in Stein's book~\cite[IV.5]{stein}.  This method is also treated in~\cite[\S 2]{ivanoplan} and in the thesis of the first author~\cite[\S7.2-3]{blairthesis}.

The multiplier theorem we use is due to Alexopolous~\cite[Theorem 6.1]{alexopolous} and treats multipliers defined with respect to the spectrum of a differential operator on a manifold, see also the work of Duong, Ouhabaz, and Sikora~\cite{duongetal}.  
It requires that the Riemannian measure is doubling and that the heat kernel $P(t,x,y)$ generated by $\Delta_\g$ should satisfy a Gaussian upper bound of the form
\begin{equation}\label{eq:heat kernel bound}
P(t,x,y) \lesssim \frac{1}{\left|B\!\left(x,\sqrt{t}\right)\right|} \, \exp\!\left(- \frac{b \, \dist_\g(x,y)^2}{t}\right)  ,
\end{equation}
where $\left|B\!\left(x,\sqrt{t}\right)\right|$ is the volume of the ball of radius $\sqrt{t}$ about $x$ and $b > 0$ is a constant.  At the end of this section, we will prove that this estimate \eqref{eq:heat kernel bound} holds on any ESCS.

Given these hypotheses, Alexopolous' theorem guarantees that any spectral multiplier $F\!\left(\sqrt{\Delta_\g}\right)$ satisfying the usual H\"ormander condition maps $L^q(X) \To L^q(X)$ for any $1 < q < \infty$.  Moreover, this boundedness is true for functions $F$ in $\mathcal{C}^N(\RR)$ which satisfy the weaker Mihlin-type condition
\begin{equation}\label{eq:Mihlin}
  \sup_{0\leq k \leq N} \sup_{\zeta \in \RR} \left|\left(\zeta \frac{d}{d\zeta}\right)^k F(\zeta)\right| \leq C < \infty  ,
\end{equation}
where $N$ is taken so that $N \geq \frac{n}{2} + 1$.

We now want to apply this theorem to a family of multipliers $F_\theta\!\left(\sqrt{\Delta_\g}\right)$, $0 \leq \theta \leq 1$, defined using the Rademacher functions $\left\{ r_m \right\}_{m=0}^\infty$.  Begin by taking
\begin{equation}
r_0(\theta) \defeq \begin{cases}
  +1, & 0 \leq \theta \leq \frac 12\\
  -1, & \frac 12 < \theta < 1  ,
\end{cases}
\end{equation}
and then extend $r_0$ to the real line by periodicity, i.e.\ $r_0(\theta + 1) = r_0 (\theta)$.  We then define the functions $r_m$ by $r_m(\theta) \defeq r_0(2^m\theta)$.  Given any square integrable sequence of scalars $\{ b_m\}_{m \geq 0}$, we consider the function $G(\theta) \defeq \sum_{m \geq 0} b_m \, r_m(\theta)$.  By a lemma in \cite[Appendix D]{stein}, for any $q$ in the interval $(1, \infty)$ there exist constants $c_q$ and $C_q$ such that
\begin{equation}
c_q \left\|G\right\|_{L^q([0,1])} \leq \left\|G\right\|_{L^2([0,1])} = \left( \sum_{m \geq 0} \left|b_m\right|^2 \right)^\frac{1}{2} \leq C_q \left\|G\right\|_{L^q([0,1])}.
\end{equation}

Define the function $\widetilde{\beta}_k(\zeta) \defeq \beta_{k-1}(\zeta)+\beta_k(\zeta)+\beta_{k+1}(\zeta)$ so that $\widetilde{\beta}_k(\zeta)\beta_k(\zeta)=\beta_k(\zeta)$. Let $F_\theta(\zeta)$ and $\widetilde{F}_\theta(\zeta)$ be the functions
\begin{equation}
F_\theta(\zeta) \defeq \sum_{k \geq 0} r_k(\theta) \, \beta_k\!\left(\sqrt{\zeta}\right) \qquad \text{and} \qquad \widetilde{F}_\theta(\zeta) \defeq \sum_{k \geq 0} r_k(\theta) \, \widetilde{\beta}_k\!\left(\sqrt{\zeta}\right)  .
\end{equation}
It can be checked that $F_\theta(\zeta)$ and $\widetilde{F}_\theta(\zeta)$ satisfy the condition~\eqref{eq:Mihlin}, and the constant $C$ appearing on the right of \eqref{eq:Mihlin} can be taken independent of $\theta$.  We thus have that for $1 < q < \infty$ and $a$ in $L^q(X)$
\begin{equation}
\begin{aligned}
\left\| \left( \sum_{k \geq 0} \left| \beta_k\!\left(\sqrt{\Delta_\g}\right) a \right|^2 \right)^{\frac{1}{2}} \right\|_{L^q(X)}^q &\lesssim \int_X \int_{\theta=0}^1 \left| \sum_{k\geq 0} r_k(\theta) \, \beta_k\!\left(\sqrt{\Delta_\g}\right) a(x) \right|^q \, \upd\theta \,\upd x \\
&\lesssim \|a\|_{L^q(X)}^q  ,
\end{aligned}
\end{equation}
and the same holds when the $\beta_k$ are replaced by the $\widetilde{\beta}_k$.

To see the other inequality in~\eqref{eq:L-P estimate}, consider $a_1$ in $L^q(X)$ and $a_2$ in $L^{q'}(X)$, and observe that
\begin{equation}
\begin{aligned}
\left|\int_X a_1 \, a_2 \, \upd x \right| &= \left|\int_X \sum_{k \geq 0} \left(\beta_k a_1\right) \left(\widetilde{\beta}_k a_2 \right) \, \upd x \right| \\
&\leq \left\|\left(\sum_{k\geq 0} \left|\beta_k \, a_1\right|^2 \right)^\frac{1}{2}\right\|_{L^q(X)} \left\|\left(\sum_{k \geq 0} \left|\widetilde{\beta}_k \, a_2 \right|^2 \right)^\frac{1}{2} \right\|_{L^{q'}(X)} \\
&\leq C \left\|\left(\sum_{k \geq 0} \left|\beta_k \, a_1\right|^2 \right)^\frac{1}{2} \right\|_{L^q(X)} \left\|a_2 \right\|_{L^{q'}(X)}  .
\end{aligned}
\end{equation}
Hence, by duality, we see that~\eqref{eq:L-P estimate} is valid.

Returning to the proof of~\eqref{eq:heat kernel bound}, we use a theorem of Grigor'yan~\cite[Theorem 1.1]{grigoryan} that establishes Gaussian upper bounds on arbitrary Riemannian manifolds.   His result implies that if $P(t,x,y)$ satisfies on-diagonal bounds
\begin{equation}\label{eq:on diagonal}
P(t,x,x) \lesssim \max\left(\frac{1}{t}, C\right)
\end{equation}
for some constant $C>0$ then there exists $b>0$ such that
\begin{equation}\label{eq:heat kernel bound2}
P(t,x,y) \lesssim \max\left(\frac{1}{t}, C\right)\, \exp\!\left(- \frac{b \, \dist_\g(x,y)^2}{t}\right)  .
\end{equation}
Since $\left|B\!\left(x,\sqrt{t}\right)\right| \approx t$ for bounded $t$, this is equivalent to~\eqref{eq:heat kernel bound}.

In order to verify~\eqref{eq:on diagonal}, we take the usual finite cover of the manifold with coordinate charts that are either isometric to a neighborhood of the Euclidean cone or the plane.  We adapt an argument of Cheeger~\cite[\S 1]{C83} to see that within each chart, the heat kernel $\widetilde{P}(t,x,y)$ of the model space is a good approximation to the intrinsic heat kernel on $X_0$.  Adjusting the cover if necessary, we may assume that the closure of any chart $Z$ is contained in a slightly larger neighborhood $Z''$ where the isometry is defined.  We may then take an intermediate neighborhood $Z'$ so that $\bar{Z} \subset Z' \subset \bar{Z'} \subset Z''$.

Let $\psi$ be a smooth cutoff supported in $Z''$ such that $\psi(z) \equiv 1$ in a neighborhood of $\bar{Z'}$.  The function $z \mapsto \psi(z)\widetilde{P}(t,z,y)$ can be seen to lie in the domain of the Laplacian $\Delta_\g$ on $X_0$, and for fixed $y$ we may consider the inhomogeneous heat equation it satisfies on that space.  For $x$ and $y$ in $Z$ we have that $\psi(x)\widetilde{P}(t,x,y) = \widetilde{P}(t,x,y)$, so Duhamel's principle shows that
\begin{equation}\label{eq:duhamel}
\widetilde{P}(t,x,y) -  P(t,x,y) =
 \int_0^t\int_{Z''} P(t-s,x,z) (\prtl_s + \Delta_\g) (\psi(z)\widetilde{P}(s,z,y))\,\upd z\,\upd s.
\end{equation}
We now recall that for fixed $y$, $\widetilde{P}(t,x,y)$ satisfies a homogeneous heat equation on the isometric space, which allows us to replace $-\prtl_s \widetilde{P}$ by the Laplacian on that space.  Applying the divergence theorem shows that~\eqref{eq:duhamel} is equal to
\begin{multline}\label{eq:gradients}
\int_0^t\int_{Z''\setminus Z'} \langle \nabla_\g P(t-s,x,z), (\nabla_\g \psi) \widetilde{P}(s,z,y) \rangle \,\upd z\, \upd s\\ - \int_0^t\int_{Z''\setminus Z'} \langle (\nabla_\g \psi) P(t-s,x,z), \nabla_\g \widetilde{P}(s,z,y)\rangle \,\upd z\, \upd s.
\end{multline}
Indeed, there is cancelation between the terms which have derivatives on both $P$ and $\widetilde{P}$.  Also, the support conditions on $\psi$ and $1-\psi$ mean that the boundary terms vanish and that the domain of integration can indeed be restricted to $Z''\setminus Z'$.

We now observe Cheeger's estimate~\cite[(1.1)]{C83}, which can be written as
\begin{equation}\label{eq:cheeger1.1}
\|d^j P(t,x,\cdot)\|_{L^2(Z''\setminus Z')} + \|d^j \widetilde{P}(t,\cdot,y)\|_{L^2(Z''\setminus Z')} \leq K_N \,t^N \qquad \text{as } t \to 0,
\end{equation}
where $j=0$ or 1 and the points $x$ and $y$ lie in $Z$.  The bound~\eqref{eq:on diagonal} now follows by proving the same estimate for $\widetilde{P}(t,x,x)$.  Indeed, for small $t$,~\eqref{eq:cheeger1.1} shows that the difference between the two kernels is negligible.  For large $t$, we observe that $Z''\setminus Z'$ can be taken to be a precompact set in either manifold.  Thus by continuity, the integrands in~\eqref{eq:gradients} can be taken to be uniformly bounded over the domain of integration and hence $P(t,x,x) \leq C$ for some large constant $C$.  When $x$ lies in a chart isometric to Euclidean space,~\eqref{eq:on diagonal} is now immediate.

To establish the on-diagonal bound for the heat kernel on the Euclidean cone, we use an approach suggested by Li~\cite[p. 284]{Li}.  A more general bound is actually announced in Theorem 2.1 of that work, but since the authors are unaware of any published proof, a weaker version of it is verified here.  In particular, we emphasize that the approach below establishes on-diagonal bounds for the Euclidean cone only.  We use $x=(r,\theta)$ to denote coordinates on the cone and to remain consistent with the notation established in~\cite{CT82},~\cite{CT82no2},~\cite{Li}, we use $\nu$ to denote the square root of the nonnegative Laplacian on the flat torus $\mathbb{R} \big/ 2\pi\rho \mathbb{Z}$.  In~\cite[(1)]{Li}, Li states the identity
\begin{multline}\label{eq:li identity}
\widetilde{P}(t,(r,\theta),(r,\theta)) =\\ \frac{1}{2\pi t} \left[ \int_0^\pi e^{-(1-\cos y)r^2/2t} \cos y\nu\;\upd y - \sin(\pi \nu) \int_0^{+\infty} e^{-(1+\cosh y)r^2/2t} e^{-y \nu}\upd y \right](\theta,\theta),
\end{multline}
where $(\theta, \theta)$ means that we are integrating the kernels of $\cos y\nu$ and $\sin \pi \nu \cdot e^{-y\nu}$ evaluated at $(\theta,\theta)$.  This identity can be verified by using Cheeger's functional calculus on cones (see e.g.~\cite[Example 3.1]{C83}) and integral representations of modified Bessel functions (see e.g. Watson~\cite[\S6.22(4)]{Watson}).

It now suffices to obtain a uniform bound on the two integrals in brackets.  For the first we use a formal identity observed by Cheeger and Taylor~\cite[(4.1), (4.8)]{CT82no2} which states that for points $\theta_1$, $\theta_2$ inside a chart on the torus, $\cos y \nu(\theta_1,\theta_2)$ is the $2\pi \rho$-periodic extension of
\begin{equation}\label{eq:point masses}
\frac 12 [\delta(\theta_1-\theta_2 + y) + \delta(\theta_1 -\theta_2 -y)].
\end{equation}
When this formal identity is made rigorous, it is subject to the proviso that if the point masses are integrated against a function with jump discontinuities, it returns the average of the left and right hand limits of the function at the center of the point mass.  Integrating $\cos y \nu(\theta,\theta)$ against the function $\mathbf{1}_{[0,\pi]}(y)\,e^{-(1-\cos y)r^2/2t}$ thus yields
\begin{equation}\label{eq:cosine eval}
\frac 12 + \sum_{k=1}^m e^{-(1-\cos y_k)r^2/2t},
\end{equation}
where $\{y_k\}_{k=1}^m$ is the (possibly empty) collection of real numbers in $(0,\pi]$ that are equivalent to 0 modulo $2 \pi \rho$.

For the second integral in~\eqref{eq:li identity}, we use the following identity in Cheeger-Taylor~\cite[(4.11)]{CT82no2} (observing that $\rho = 1/\gamma$)
\begin{equation}\label{eq:sine exp}
\sin \nu \pi \cdot e^{-\nu y}(\theta,\theta)
=\frac 1{2\pi \rho} \left( \frac{\sin(\pi/\rho)}{ \cosh (y/\rho) - \cos ( \pi/\rho)}\right).
\end{equation}
This gives rise to the integral
\begin{equation}\label{eq:sine exp integral}
\frac 1{2\pi \rho} \int_0^\infty e^{-(1+\cosh y)r^2/2t} \frac{\sin(\pi/\rho)}{ \cosh (y/\rho) - \cos ( \pi/\rho)}\;dy.
\end{equation}
Note that this integral vanishes when $\rho = 1/N$ for $N$ a positive
integer; this corresponds to the absence of diffraction on cones of
these radii.  Otherwise, the integrand is bounded near 0 and rapidly decaying at infinity.  This provides uniform bounds on the second integral.


\end{document}